\documentclass[final]{siamltex}
\newcommand{\Gn}{{\cal G}_n}
\newcommand{\Aut}{{\it Aut}}
\newcommand{\Prob}{{\it Prob}}

\newcommand{\Avg}{{\it Avg}}
\newcommand{\In}{\{0,1,\ldots,n-1\}}
\newtheorem{claim}{\sc Claim}
\title{Kolmogorov Random Graphs and the Incompressibility Method\thanks{A
partial preliminary version appeared in the
{\em Proc. Conference on Compression and Complexity of Sequences},
IEEE Comp. Sci. Press, 1997. HB, JT and PV were 
partially supported by the European Union
through NeuroCOLT ESPRIT Working Group Nr. 8556,
and by  NWO through NFI Project ALADDIN
number NF 62-376 and 
SION Grant 612-34-002;
ML was supported in part by
the NSERC Operating Grant OGP0046506, CITO, a CGAT grant, and the
Steacie Fellowship.
}}

\author{Harry Buhrman\thanks{
CWI,
Kruislaan 413, 1098 SJ Amsterdam, The Netherlands.
Email: buhrman@cwi.nl}
\and
Ming Li\thanks{
Department of Computer Science, University of Waterloo,
Waterloo, Ont. N2L 3G1, Canada. E-mail: mli@wh.math.uwaterloo.ca}
\and
John Tromp\thanks{
CWI,
Kruislaan 413, 1098 SJ Amsterdam, The Netherlands.
Email: tromp@cwi.nl}
\and
Paul Vit\'{a}nyi\thanks{
CWI,
Kruislaan 413, 1098 SJ Amsterdam, The Netherlands.
Email: paulv@cwi.nl}
}
\date{}

\begin{document}

\maketitle
 
\begin{abstract}
We investigate topological, combinatorial, statistical, and enumeration
properties of finite graphs with high Kolmogorov complexity
(almost all graphs) 
using the novel incompressibility method. Example results are:
(i) the mean and variance of the number 
of (possibly overlapping) ordered labeled subgraphs of a labeled graph
as a function of its randomness deficiency (how far it falls short of
the maximum possible Kolmogorov complexity)
and (ii) a new elementary proof for the number of unlabeled graphs.
\end{abstract}



\pagestyle{myheadings}
\thispagestyle{plain}
\markboth{H. BUHRMAN, M. LI, J. TROMP AND P. VITANYI}{KOLMOGOROV RANDOM GRAPHS}

\section{Introduction}
The incompressibility of individual random objects  
yields a simple but powerful proof technique.
The incompressibility method,
\cite{LiVi93}, is a new general purpose
tool and should be compared with
the pigeon hole principle\index{pigeon hole principle} or the
probabilistic method\index{probabilistic method}. 
Here we apply the incompressibility method to randomly generated graphs
and ``individually random'' graphs---graphs with high Kolmogorov complexity.

In a typical proof using the incompressibility method,
one first chooses an individually random object from the
class under discussion.
This object is
effectively incompressible.
The argument invariably says that if a desired property
does not hold, then the object
can be compressed. This yields the required contradiction.
Since
a randomly generated object is 
{\em with overwhelming probability} 
individually random and hence incompressible,
one usually obtains the property with high probability.

{\bf Results}
We apply the incompressibility
method to obtain combinatorial properties of graphs with high Kolmogorov
complexity. These properties are parametri\-zed in terms
of a ``randomness deficiency'' function.\footnote{Randomness deficiency
measures how far the object falls short of the maximum possible
Kolmogorov complexity. It is formally defined in Definition~\ref{def.rg}.}
 This can be considered
as a parametri\-zed version of the incompressibility method.  In
Section~\ref{sect.topol} we show that: 
For every labeled graph on $n$ nodes with high Kolmogorov
complexity (also called ``Kolmogorov random 
graph'' or ``high complexity graph'')
the node degree of every vertex is about $n/2$ and 
there are about $n/4$ node-disjoint paths of length 2 
between every pair of nodes.
In Section~\ref{sect.statsubgr}, we analyze `normality'
properties of Kolmogorov random graphs. In analogy with infinite
sequences one can call an infinite labeled graph `normal' if each
finite ordered labeled subgraph of size $k$ occurs in the appropriate sense
(possibly overlapping) with
limiting frequency $2^{-{k \choose 2}}$. It follows from
Martin-L\"of's theory of effective tests for randomness \cite{Ma66}
that individually random (high complexity) infinite labeled graphs are
normal.
Such properties cannot hold precisely for finite graphs, where 
randomness is necessarily a matter of degree: We determine close
quantitative bounds on the normality (frequency of subgraphs) of high complexity
finite graphs in terms of their randomness deficiency.
Denote the number of
unlabeled graphs on $n$ nodes by $g_n$.
In Section~\ref{sect.unlabeled} we demonstrate the use
of the incompressibility method and Kolmogorov random graphs by
providing a new elementary proof that $g_n \sim  2^{n \choose 2} / n!$.
This has previously been obtained
by more advanced methods, \cite{HP73}. Moreover, we give a good
estimate of the error term.
Part of the proof involves estimating the order (number of automorphisms)
$s(G)$ of graphs $G$ as a function of the randomness deficiency
of $G$. For example, we show that labeled
graphs with randomness deficiency appropriately less than $n$
are rigid (have but one automorphism: the identity automorphism).

{\bf Related Work}
Several properties above (high degree nodes, diameter 2,
rigidity) have also been proven by traditional
methods to hold with high probability
for randomly generated graphs, \cite{Bo85}. We provide new proofs
for these results
using the incompressibility method. They are actually
proved to hold for the definite class of
Kolmogorov random graphs---rather than with high probability
for randomly generated graphs. 
In \cite{LiVi94b}
(also \cite{LiVi93}) two of us investigated
topological properties of 
labeled graphs with high Kolmogorov complexity
and proved them using the incompressibility
method to compare ease of such proofs with the probabilistic
method \cite{ES74} and entropy method.
In \cite{Ki92}
it was shown that every labeled tree on $n$ nodes with randomness
deficiency $O(\log n)$ has
maximum node degree of $O( \log n / \log \log n)$.
Analysis of Kolmogorov random graphs was used to establish
the total interconnect length
of Euclidean (real-world) embeddings of computer network
topologies \cite{Vi95}, 
and the size of compact routing tables in computer networks
\cite{BHV95}.
Infinite binary sequences that asymptotically have
equal numbers of 0's and 1's, and more generally, 
where every block of length $k$
occurs (possibly overlapping) with frequency 
$1/2^k$ 
were called  ``normal'' by E. Borel, \cite{Bo14}. 
References \cite{LiVi93,LiVi94a} 
investigate 
the quantitative deviation from normal 
as a function of the Kolmogorov complexity of a finite binary string.
Here we consider analogous question 
for Kolmogorov random graphs.
\footnote{There are some results 
along these lines related to randomly generated graphs, 
but as far as the
authors could
ascertain
(consulting Alan Frieze, Svante Janson, and
Andrzej Rucinski around June 1996)
such properties 
have not been investigated in the same detail as here.
See for example
\cite{ASE92}, pp. 125--140. But note that 
also pseudo-randomness is different
from Kolmogorov randomness.
}
Finally, there is a close relation and
genuine differences between 
high-probability properties
and properties of incompressible objects,
see  \cite{LiVi93}, Section 6.2.

\subsection{Kolmogorov complexity}
We use the following notation.
Let $A$ be a finite set. By $d(A)$ we denote the {\em cardinality}
of $A$. In particular, $d(\emptyset)=0$. Let $x$ be
a finite binary string. Then $l(x)$ denotes the {\em length}
(number of bits) of $x$. In particular, $l(\epsilon)=0$
where $\epsilon$ denotes the {\em empty word}.

Let $x,y,z \in {\cal N}$, where
${\cal N}$ denotes the natural
numbers. Identify
${\cal N}$ and $\{0,1\}^*$ according to the
correspondence 
\[(0, \epsilon ), (1,0), (2,1), (3,00), (4,01), \ldots . \]
Hence, the length $l(x)$ of $x$ is the number of bits
in the binary string or number $x$.
Let $T_0 ,T_1 , \ldots$ be a standard enumeration
of all Turing machines.
Let $\langle \cdot ,\cdot \rangle$ be a standard one-one mapping
from ${\cal N} \times {\cal N}$
to ${\cal N}$, for technical reasons choosen such that
$l(\langle x ,y \rangle) = l(y)+O(l(x))$.
An example is $\langle x ,y \rangle = 1^{l(x)}0xy$.
This can be iterated to
$\langle  \langle \cdot , \cdot \rangle , \cdot \rangle$.

Informally, the Kolmogorov complexity, \cite{Ko65},
of $x$ is the length of the
{\em shortest} effective description of $x$.
That is, the {\em Kolmogorov complexity} $C(x)$ of
a finite string $x$ is simply the length
of the shortest program, say in
FORTRAN (or in Turing machine codes)
encoded in binary, which prints $x$ without any input.
A similar definition holds conditionally, in the sense that
$C(x|y)$ is the length of the shortest binary program
which computes $x$ on input $y$. 
Kolmogorov complexity is absolute in the sense
of being independent of the programming language,
up to a fixed additional constant term which depends on the programming
language but not on $x$. We now fix one canonical programming
language once and for all as reference and thereby $C()$.
For the theory and applications, see \cite{LiVi93}.
A formal definition is as follows:

\begin{definition}
\rm
Let $U$ be an appropriate universal Turing machine
such that 
\[U(\langle \langle i,p \rangle ,y \rangle ) =
T_i (\langle p,y\rangle) \]
 for all $i$ and $\langle p,y\rangle$.
The {\em conditional Kolmogorov complexity} of $x$ given $y$
is
\[C(x|y) = \min_{p \in \{0,1\}^*} \{l(p): U (\langle p,y\rangle)=x \}. \]
The unconditional Kolmogorov complexity of $x$ is defined
as $C(x) := C(x| \epsilon )$.
\end{definition}
It is easy to see that there are strings that can be described
by programs much shorter than themselves. For instance, the
function defined by $f(1) = 2$ and $f(i) = 2^{f(i-1)}$
for $i>1$ grows very fast, $f(k)$ is a ``stack'' of $k$ twos.
Yet for each $k$ it is clear that $f(k)$
has complexity at most $C(k) + O(1)$.
What about incompressibility?

By a simple counting argument one can show
that whereas some strings can be enormously compressed,
the majority of strings can hardly be compressed
at all.
 
For each $n$ there are $2^n$ binary
strings of length $n$, but only
$\sum_{i=0}^{n-1} 2^i = 2^n -1$ possible shorter descriptions.
Therefore, there is at least one binary string
$x$ of length $n$ such that $C(x)   \geq   n$.
We call such strings $incompressible$. It also
follows that for any length $n$ and any binary string $y$,
there is a binary string $x$ of length $n$ such that
$C(x| y)   \geq   n$. Generally,
for every constant $c$ we can call a string $x$ is
\it c-incompressible 
\rm if $C(x)   \geq   l(x) -c$.
Strings that are incompressible (say, $c$-incompressible
with small $c$) are patternless,
since a pattern could be used to reduce
the description length. Intuitively, we
think of
such patternless sequences
as being random, and we
use ``random sequence''
synonymously with ``incompressible sequence.''
\footnote{It is possible to give a rigorous 
formalization of the intuitive notion
of a random sequence as a sequence that passes all
effective tests for randomness, see for example \cite{LiVi93}.
}
By the same counting argument as before we find that the number
of strings of length $n$ that are $c$-incompressible
is at least $2^n - 2^{n-c} +1$. Hence
there is at least one 0-incompressible string of length $n$,
at least one-half of all strings of length $n$ are 1-incompressible,
at least three-fourths  of all strings
of length $n$ are 2-incompressible, \ldots , and
at least the $(1- 1/2^c )$th part
of all $2^n$ strings of length $n$ are $c$-incompressible. This means
that for each constant $c   \geq   1$ the majority of all
strings of length $n$ (with $n   >   c$) is $c$-incompressible.
We generalize this to the following simple but extremely
useful:
\begin{lemma}
\label{C2}
Let $c$ be a positive integer.
For each fixed $y$, every
set $A$ of cardinality $m$ has at least $m(1 - 2^{-c} ) + 1$
elements $x$ with $C(x| y)   \geq   \lfloor \log m \rfloor  - c$.
\end{lemma}
\begin{proof}
By simple counting.
\end{proof}
 
As an example, set $A =  \{ x: l(x) = n  \}  $. Then the cardinality
of $A$ is $m = 2^n$.
Since it is easy to assert that $C(x)  \leq n + c$ for some
fixed $c$ and all $x$ in $A$, Lemma~\ref{C2} demonstrates
that this trivial estimate is quite sharp. The deeper
reason is that since there are few short programs, there can
be only few objects of low complexity. We require another quantity:
the prefix Kolmogorov complexity which is defined just as $C( \cdot | \cdot)$
but now with respect to a subset of Turing machines that have
the property that the set of programs for which the machine halts
is prefix-free, that is, no halting program is a prefix of any other
halting program. For details see \cite{LiVi93}. Here we require
only the quantitative relation below.
\begin{definition}
The {\em prefix} Kolmogorov complexity of $x$ conditional to $y$
is denoted by $K(x|y)$. It satisfies the inequality
\[C(x|y) \leq K(x|y) \leq C(x|y) + 2 \log C(x|y) + O(1) .\]
\end{definition}

\section{Kolmogorov Random Graphs}\label{sect.KRgraphs}
\label{sect.topol}
Statistical properties of strings with high Kolmogorov complexity
have been studied in \cite{LiVi94a}.
The interpretation of strings as
more complex combinatorial objects leads to a 
new set of properties and problems that have no direct
counterpart in the ``flatter'' string world. Here we 
derive topological, combinatorial, and statistical
properties of graphs with high Kolmogorov complexity.
Every such graph
possesses simultaneously all properties 
that hold with high probability for randomly generated graphs.
They constitute ``almost all graphs'' and the derived properties
a fortiori hold with probability that goes to 1
as the number of nodes grows unboundedly.
\begin{definition}\label{def.gc}
\rm
Each labeled graph $G=(V,E)$ on $n$ nodes $V=\{1,2,\ldots, n\}$
can be represented (up to automorphism)
by a binary string $E(G)$ of length ${n \choose 2}$. We
simply assume a fixed ordering of the
${n \choose 2}$ possible edges in an $n$-node graph, e.g.
lexicographically, and
let the $i$th bit in the string indicate presence (1) or absence (0) of
the $i$'th edge. Conversely, each
binary string of length ${n \choose 2}$ encodes an $n$-node graph.
Hence we can identify each such graph with
its binary string representation.
\end{definition}
\begin{definition}\label{def.rg}
\rm
A labeled graph $G$ on $n$ 
nodes has {\em randomness deficiency}\index{randomness deficiency}
at most $\delta (n)$, and is called
$\delta (n)$-{\em random},
if it satisfies
\begin{equation}\label{eq.KG}
 C(E(G)|n ) \geq {n \choose 2} - \delta (n).
\end{equation}
\end{definition}

\subsection{Some Basic Properties}
Using Lemma~\ref{C2}, 
with $y=n$, $A$ the set of strings of length $n \choose 2$, and
$c=\delta(n)$ gives us
\begin{lemma}\label{lem.frac}
A fraction of at least
\label{eq.count}
$1 - 1/2^{\delta (n)}$
of all labeled  graphs $G$ on $n$ nodes is $\delta (n)$-random.
\end{lemma}
As a consequence, for example the $c \log n$-random
labeled graphs constitute a fraction of
at least $(1 - 1/n^c)$ of all graphs on
$n$ nodes, where $c>0$ is an arbitrary 
constant.

Labeled graphs with high-complexity  have many specific topological
properties, which seem to contradict their randomness.
However, these are simply the likely properties, whose absence would
be rather unlikely.
Thus, randomness enforces strict statistical regularities.
For example, to have diameter exactly two.

We will use the following lemma (Theorem 2.6.1 in \cite{LiVi93}):
\begin{lemma}\label{blockszerotex}
Let $x=x_1 \ldots x_n$ be a binary string of
length $n$, and $y$ a much smaller string
of length $l$. Let $p = 2^{-l}$ and
$\#y(x)$ be the number of
(possibly overlapping) distinct occurrences of $y$ in $x$.
For convenience, we assume that
$x$ ``wraps around'' so that an occurrence
of $y$ starting at the end of $x$
and continuing at the start also counts.
Assume that $l \leq \log n$.
There is a constant $c$
such that for all $n$ and $x \in \{0,1\}^n$,
if $C(x) \geq n - \delta(n)$, then
\[ |\#y(x)-pn|  \leq  \sqrt{ \alpha pn},\]
with $\alpha = [K(y|n)+\log l + \delta(n)+c] 3l / \log e $.
\end{lemma}

\begin{lemma}\label{lem.diam}
All $o(n)$-random labeled graphs have $n/4+o(n)$
disjoint paths of length 2 between each pair of nodes $i,j$.
In particular, all $o(n)$-random  labeled graphs
have diameter 2.
\end{lemma}
\begin{proof}
The only graphs with diameter 1 are the complete graphs that
can be described in $O(1)$ bits, given $n$, and hence are not random.
It remains to consider an  $o(n)$-random graph $G=(V,E)$ with diameter
greater than or equal to 2. Let $i,j$  be a pair of nodes connected 
by $r$ disjoint paths of length 2.
Then we can describe $G$ by modifying the old
code for $G$ as follows:
\begin{itemize}
\item
A program to reconstruct the object from the various
parts of the encoding in $O(1)$ bits;
\item
The identities of $i < j$ in  $2  \log  n$ bits;
\item
The old code $E(G)$ of $G$ with the $2(n-2)$ bits representing
presence or absence of edges $(j,k)$ and $(i,k)$ for each 
$k \neq i,j$ deleted.
\item a shortest program for the string $e_{i,j}$ consisting of the (reordered)
$n-2$ pairs of bits deleted above.
\end{itemize}
>From this description we can reconstruct $G$ in
\[ O(\log n) + {n \choose 2}  - 2(n-2) + C(e_{i,j}|n)  \]
bits, from which we may conclude that $C(e_{i,j}|n) \geq l(e_{i,j}) - o(n)$.
As shown in \cite{LiVi94a} or \cite{LiVi93} (here Lemma~\ref{blockszerotex})
this implies that the frequency of occurrence in $e_{i,j}$
of the aligned 2-bit block `11'---which by construction equals the number of
disjoint paths of length 2 between $i$ and $j$---is $n/4 + o(n)$.

\end{proof}

A graph is {\em $k$-connected} if there are at least $k$ node-disjoint
paths between every pair of nodes.

\begin{corollary}
All $o(n)$-random labeled graphs are $( \frac{n}{4}+o(n))$-connected.
\end{corollary}

\begin{lemma}
Let $G=(V,E)$ be a graph on $n$ nodes with randomness deficiency
$O(\log n)$. Then the largest clique in $G$ has at most
$\lfloor 2 \log n \rfloor + O(1)$ nodes.
\end{lemma}
\begin{proof}
Same proof as largest size transitive subtournament in 
high complexity tournament as in \cite{LiVi93}.
\end{proof}

With respect to the related property of random graphs, in \cite{ASE92}
pp. 86,87  it is shown that a random graph with edge probability
$1/2$ contains a clique on asymptotically $2 \log n$ nodes with probability
at least $1-e^{-n^2}$.

\subsection{Statistics of Subgraphs}\label{sect.statsubgr}
We start by defining the notion of labeled subgraph of a labeled graph.
\begin{definition}
\rm
Let $G=(V,E)$ be a labeled graph on $n$ nodes.
Consider a labeled graph 
$H$ on $k$ nodes $\{1,2, \ldots,k\}$. 
Each 
subset of $k$ nodes of $G$ induces a 
subgraph $G_k$ of $G$. The subgraph $G_k$
is an  ordered labeled {\em occurrence} of $H$  when we obtain $H$ by
relabeling the nodes $i_1< i_2<  \cdots < i_k$ of $G_k$ as
$1,2, \ldots, k$. 
\end{definition}

It is easy to conclude from the statistics of high-complexity
strings in 
Lemma~\ref{blockszerotex}
that the frequency of each of the two labeled two-node subgraphs
(there are only two different ones: the graph consisting
of two isolated nodes and the graph consisting of 
two connected nodes) 
in a $\delta (n)$-random graph $G$ is 
\[ \frac{ n(n-1)}{4} \pm \sqrt{ \frac{3}{4}(\delta (n)+ O(1)) n(n-1)/ \log e}. \]
This case is easy since the frequency of such subgraphs
corresponds to the frequency of 1's or $0$'s in the ${n \choose 2}$-length
standard encoding $E(G)$ of $G$. 
However, to
determine the frequencies of labeled subgraphs
on $k$ nodes (up to isomorphism) for $k>2$ is a matter more complicated 
than the frequencies of substrings of length $k$.
Clearly, there are $n \choose k$ subsets of $k$ nodes out of $n$
and hence that many occurrences of subgraphs. Such subgraphs may
overlap in more complex ways than substrings of a string. 
Let $\#H(G)$ be {\em the number of times $H$ occurs}
as an ordered labeled subgraph of $G$ (possibly overlapping). Let 
$p$ be the probability that we obtain $H$ by flipping a fair
coin to decide for each pair of nodes whether
it is connected by an edge or not,
\begin{equation}\label{eq.defp}
p=2^{-k(k-1)/2}.
\end{equation}
\begin{theorem}\label{theo.freqG}
Assume the terminology above with 
$G=(V,E)$ a labeled  graph on $n$ nodes, $k$ is a positive integer
dividing $n$, and
$H$ is a labeled graph on $k \leq \sqrt{2 \log n}$ nodes. Let
$C(E(G)|n ) \geq {n \choose 2} - \delta (n)$. Then
\[ \left|\#H(G)- {n \choose k}p \right| \leq 
 {n \choose k} \sqrt{\alpha (k/n) p} , \]
with
$\alpha := (K(H|n) + \delta(n) + \log {n \choose k}/(n/k) + O(1))
3 / \log e$. 
\end{theorem}
\begin{proof}
A {\em cover} of $G$ is a
set $C= \{S_1, \ldots , S_N\}$ with $N=n/k$, where the $S_i$'s are pairwise
disjoint subsets of  $V$
and $\bigcup_{i=1}^N S_i = V$.
According to  \cite{Ba74}:
\begin{claim}\label{claim.bara}
\rm
There is a partition of the ${n \choose k}$ different $k$-node subsets
into $h={n \choose k}/N$ distinct covers of $G$, 
each cover consisting of $N = n/k$ disjoint subsets. That is,
each subset of $k$ nodes of $V$ belongs to precisely one cover.
\end{claim}
Enumerate the covers as $C_0, C_2 , \ldots , C_{h-1}$.
For each $i \in \{ 0,1, \ldots, h-1 \}$ 
and $k$-node labeled graph $H$, let $\#H(G,i)$ be the number
of (now non overlapping) occurrences of subgraph $H$ in $G$
occurring in cover $C_i$.

Now consider an experiment of $N$ trials, each trial with
the same set of $2^{k(k-1)/2}$ 
outcomes. 
Intuitively, each trial
corresponds to an element of a cover, and each outcome
corresponds to a $k$-node subgraph.
For every $i$ we can form a string $s_i$ consisting of the
$N$ blocks of ${k \choose 2}$ bits that represent presence
or absence of edges within the induced subgraphs 
of each of the $N$ subsets of $C_i$.
Since $G$ can be reconstructed from $n,i,s_i$ and the remaining 
${n \choose 2} - N {k \choose 2}$ bits of $E(G)$,
we find that $C(s_i|n) \geq l(s_i) - \delta (n) - \log h$.
Again, according to Lemma~\ref{blockszerotex} this implies
that the frequency of occurrence of the aligned ${k \choose 2}$-block
$E(H)$, which is $\#H(G,i)$, equals
\[ Np \pm  \sqrt{Np \alpha},  
\]
with $\alpha$ as in the theorem statement.
One can do this for each $i$ independently,
notwithstanding the dependence between
the frequencies of subgraphs in different covers. Namely, the argument
depends on the incompressibility of $G$ alone. If the number
of occurrences of a certain subgraph in {\em any} of the
covers is too large or too small then we can compress $G$.
Now,
\begin{eqnarray*}
 \left| \#H(G) - p{n \choose k} \right| & = &
\sum_{i=0}^{h-1} |\#H(G,i)-Np| \\
& \leq & {n \choose k} \sqrt{\alpha (k/n)p}.
\end{eqnarray*}
\end{proof}

In \cite{LiVi93,LiVi94a} we investigated up to which
length $l$ all blocks of length $l$ occurred at least once in
each $\delta (n)$-random string of length $n$.
\begin{theorem}
Let $\delta (n) <  2^{\sqrt{\frac{1}{2} \log n}}/4 \log n $
and  $G$ be a $\delta (n)$-random graph on $n$ nodes.
Then for sufficiently large $n$, the graph $G$
contains all subgraphs on
$\sqrt{ 2 \log n}$  nodes.
\end{theorem}
\begin{proof}
We are sure that $H$ on $k$ nodes occurs at least once in $G$
if $ {n \choose k} \sqrt{ \alpha (k/n) p}$
in Theorem~\ref{theo.freqG} is less than ${n \choose k}p$.
This is the case if $\alpha < (n/ k) p$.
This inequality is satisfied for an overestimate of
$K(H|n)$ by 
  ${k \choose 2} + 2 \log {k \choose 2} +O(1)$
(since $K(H|n) \leq K(H)+O(1)$), and $p=2^{-k(k-1)/2}$,
and $k$ set at
$k = \sqrt{ 2 \log n}$. This proves the theorem.
\end{proof}
\index{statistical properties!of graphs|)}

\subsection{Unlabeled Graph Counting}
\label{sect.unlabeled}
An unlabeled graph is a graph with no labels. For convenience we can
define this as follows: Call two labeled graphs {\em equivalent}
(up to relabeling) if there is a relabeling that makes them equal.
An {\em unlabeled graph} is an equivalence class of labeled graphs.
An {\em automorphism} of $G=(V,E)$ is a permutation $\pi$ of $V$
such that $(\pi(u),\pi(v)) \in E$ iff $(u,v)\in E$.
Clearly, the set of automorphisms of a graph forms a group
with group operation of function composition and
the identity permutation as unity. It is easy to verify that
$\pi$ is an automorphism of $G$  iff $\pi (G)$ and $G$
have the {\em same binary string standard encoding}, that is,
$E(G)=E(\pi (G))$. This contrasts
with the more general case of permutation relabeling, where the
standard encodings may be different. 
A graph is {\em rigid} if its only automorphism is the identity automorphism.
It turns out that Kolmogorov random graphs are rigid graphs.
To obtain an expression for the number
of unlabeled graphs we have to estimate the number
of automorphisms of a graph in terms of its randomness deficiency.

In \cite{HP73} an asymptotic expression for the number of unlabeled 
graphs is derived using sophisticated methods.
We give a new elementary proof by incompressibility.
Denote by $g_n$ the number of unlabeled graphs on $n$
nodes---that is, the number of isomorphism classes in the set $\Gn$
of undirected graphs on nodes $\In$.

\begin{theorem}\label{theo.unlabeled}
$g_n \sim \frac{2^{n \choose 2}}{n!}$. 
\end{theorem}

\begin{proof}
Clearly, 
\[ g_n = \sum_{G \in \Gn} \frac{1}{d( \bar{G})}, \]
 where $\bar{G}$ is the
isomorphism class of graph $G$. By elementary group theory,
\[ d(\bar{G}) = \frac{d(S_n)}{d(\Aut(G))} = \frac{n!}{d(\Aut(G))}, \]
where $S_n$ is the group of permutations on $n$ elements, and $\Aut(G)$
is the automorphism group of $G$.
Let us partition $\Gn$ into $\Gn = \Gn^0 \cup \ldots \cup \Gn^n$,
where $\Gn^m$ is the set of graphs for which $m$ is the number of nodes
moved (mapped to another node) by any of its automorphisms.

\begin{claim}
For $G \in \Gn^m$, $d( \Aut(G)) \leq n^m = 2^{m\log n}$.
\end{claim}

\begin{proof}
$d( \Aut(G)) \leq {n \choose m}m! \leq n^m$.
\end{proof}

Consider each graph $G \in \Gn$ having a probability
$\Prob(G) = 2^{-{n \choose 2}}$.

\begin{claim}
$\Prob(G \in \Gn^m) \leq 2^{-m(\frac{n}{2}-\frac{3m}{8} -\log n)}$.
\end{claim}

\begin{proof}
By Lemma~\ref{lem.frac} it suffices to show that, if
$G \in \Gn^m$ and
\[ C(E(G)|n,m) \geq {n \choose 2} - \delta (n,m) \]
then $\delta (n,m)$ 
satisfies
\begin{equation}\label{eq.dnm}
\delta (n,m) \geq m(\frac{n}{2}-\frac{3m}{8} -\log n). 
\end{equation}
Let $\pi \in \Aut(G)$ move $m$ nodes. Suppose $\pi$ is the product of $k$
disjoint cycles of sizes $c_1,\ldots, c_k$.
Spend at most $m \log n$ bits describing $\pi$: 
For example, if the nodes $i_1 < \cdots < i_m$ are moved 
then list the sequence $\pi(i_1),\ldots, \pi(i_m)$. Writing
the nodes of the latter sequence in increasing order 
we obtain $i_1 , \dots , i_m$ again, that is, we execute permutation $\pi^{-1}$,
and hence we obtain $\pi$.

Select one node from each cycle---say, the lowest numbered one.
Then for every unselected node on a cycle, 
we can delete the $n-m$ bits corresponding
to the presence or absence of edges to stable nodes, 
and $m-k$ half-bits corresponding to presence or absence of edges to the other,
unselected cycle nodes.
In total we delete
\[ \sum_{i=1}^{k} (c_i -1)(n-m + \frac{m-k}{2}) = (m-k)(n- \frac{m+k}{2}) \]
bits. Observing that $k=m/2$ is the largest possible value for $k$,
we arrive at the claimed 
$\delta (n,m)$ of $G$ (difference between savings and spendings is
$\frac{m}{2}(n- \frac{3m}{4}) - m \log n$)
of Equation~\ref{eq.dnm}.
\end{proof}

We continue the proof of the main theorem:
\[g_n = \sum_{G \in \Gn} \frac{1}{d(\bar{G})}
= \sum_{G \in \Gn} \frac{d(\Aut(g))}{n!}
= \frac{2^{{n \choose 2}}}{n!} E_n,
\]
where
$E_n := \sum_{G \in \Gn} \Prob(G) d(\Aut(G))$
is the expected size of the automorphism group of a graph on $n$ nodes.
Clearly, $E_n \geq 1$, yielding the lower bound on $g_n$.
For the upper bound on $g_n$, noting that $\Gn^1 = \emptyset$ and using
the above claims, we find
\begin{eqnarray*}
E_n & = & \sum_{m=0}^n \Prob(G \in \Gn^m) \Avg_{G \in \Gn^m} d(\Aut(G))  \\
& \leq & 1 + \sum_{m=2}^n 2^{-m(\frac{n}{2}-\frac{3m}{8}  -2 \log n)} \\
& \leq & 1 + 2^{-(n - 4 \log n - 2)}.
\end{eqnarray*}
which proves the theorem.
\end{proof}

The proof of the theorem shows that the error in the asymptotic
expression is very small:
\begin{corollary}
$\frac{2^{{n \choose 2}}}{n!} \leq g_n \leq 
\frac{2^{{n \choose 2}}}{n!} (1+ \frac{4n^4}{2^n})$.
\end{corollary}

It follows from Equation~\ref{eq.dnm} that (since $m=1$ is impossible):
\begin{corollary}
If a graph $G$ has randomness deficiency slightly less than $n$
(more precisely, $C(E(G)|n) \geq {n \choose 2} - n - \log n -2$)
then $G$ is rigid.
\end{corollary}

The expression for $g_n$ can be used to determine
the maximal complexity of an unlabeled graph on $n$ nodes.
Namely, we can effectively enumerate all unlabeled graphs as follows:
\begin{itemize}
\item
Effectively enumerate all labeled graphs on $n$ nodes
by enumerating all binary strings of length $n$ and
for each labeled graph $G$ do:
\subitem
If $G$ cannot be obtained by relabeling from any previously 
enumerated labeled graph then $G$ is added to the set of
unlabeled graphs.
\end{itemize}
This way we obtain each unlabeled graph
by precisely one labeled graph representing it.
Since we can describe each unlabeled graph by its index
in this enumeration, we find
by Theorem~\ref{theo.unlabeled} and 
Stirling's formula\index{Stirling's formula} that
if $G$ is an unlabeled graph then
\[ C(E(G)|n) \leq {n \choose 2} - n \log n + O(n) . \]
\begin{theorem}\label{theo.drop}
Let $G$ be a labeled graph on $n$ nodes and let $G_0$ be
the unlabeled version of $G$. There exists a graph
$G'$ and a label permutation $\pi$ such that $G' = \pi(G)$
and up to additional constant terms $C(E(G'))=C(E(G_0))$ and
$C(E(G)|n) =  C(E(G_0) , \pi |n )$.
\end{theorem}

By Theorem~\ref{theo.drop},
for {\em every} graph $G$ on $n$ nodes with maxi\-mum-complexity 
there is a relabeling (permutation)
that causes the complexity to drop by as much as
$n \log n$. Our proofs of topological properties
 by the incompressibility method 
required the graph $G$ to be Kolmogorov random in the sense of
$C(E(G)|n) \geq {n \choose 2} - O(\log n)$ or for some results
$C(E(G)|n) \geq {n \choose 2} - o(n)$. 
Hence by relabeling such a graph
we can always obtain a labeled graph that 
has a complexity too low to use our incompressibility proof.
Nonetheless, topological properties do not change under
relabeling.


\begin{thebibliography}{99}

\bibitem{ASE92}
{\sc N. Alon, J.H. Spencer
and P. Erd\H{o}s}, The Probabilistic Method,
Wiley, 1992, 

\bibitem{Bo14}
{\sc E. Borel}, 
Le\c{c}ons sur la th\'eorie des functions, 2nd Edition,
1914, 182--216.

\bibitem{Ba74}
{\sc Baranyai, Zs.}, 
{\em On the factorization of the complete uniform hypergraph,}
pp. 91-108 in: A. Hajnal, R. Rado, V.T. S\'os, Eds, 
Infinite and Finite Sets,  Proc. Coll. Keszthely,
Colloq. Math. Soc. J\'{a}nos Bolyai, Vol. 10 (1995)
North-Holland, Amsterdam.

\bibitem{Bo79}
{\sc B. Bollob\'as}, Graph Theory,
Springer-Verlag, New York, 1979.

\bibitem{Bo85}
{\sc B. Bollob\'as}, Random Graphs,
Academic Press, London, 1985.

\bibitem{BHV95}
{\sc H. Buhrman, J.H. Hoepman, and P. Vit\'anyi},
{\em Space-Efficient Routing Tables for Almost All Networks
 and the Incompressibility Method}, {\em SIAM J. Comput.},
To appear. 

\bibitem{ES74}
{\sc P. Erd\"os and J. Spencer}, Probabilistic Methods in Combinatorics,
Academic Press, New York, 1974.
\bibitem{Ki92}
{\sc W.W. Kirchherr},
{\em Kolmogorov complexity and random graphs},
Inform.\ Process. Lett., 41(1992), 125--130.
\bibitem{LiVi93}
{\sc M. Li and P.M.B. Vit\'anyi},
An Introduction to Kolmogorov Complexity
and its Applications,
Springer-Verlag, New York, 2nd Edition, 1997.

\bibitem{LiVi94b}
{\sc M. Li and P.M.B. Vit\'anyi}, 
{\em Kolmogorov complexity arguments in
Combinatorics}, J. Comb. Th., Series A, 66:2(1994), 226-236.
Errata, {\em Ibid.}, 69(1995), 183.

\bibitem{LiVi94a}
{\sc M. Li and P.M.B. Vit\'anyi}, {\em Statistical properties of finite sequences
with high Kolmogorov complexity},
Math. System Theory, 27(1994), 365-376.

\bibitem{HP73}
{\sc F. Harary and E.M. Palmer},
Graphical Enumeration,
Academic Press, 1973.

\bibitem{Ko65}
{\sc A.N. Kolmogorov},
{\em Three approaches to the quantitative definition of information}.
Problems Inform. Transmission, 1(1):1--7, 1965.

\bibitem{Ma66}
{\sc P. Martin-L\"{o}f}, {\em On the definition of random sequences},
Information and Control, 9(1966), 602-619.

\bibitem{Vi95}
{\sc P.M.B. Vit\'anyi}, {\em Physics and the New Computation},
Prague, August 1995,
Proc. 20th Int. Symp. Math. Foundations of Computer Science, MFCS'95,
Lecture Notes in Computer Science, Vol 969, Springer-Verlag,
Heidelberg, 1995,  106--128.
\end{thebibliography}
\end{document}